\documentclass{amsart}
\usepackage{amsmath}%
\usepackage{amsfonts}%
\usepackage{amssymb}%
\usepackage{graphicx}
\usepackage{times}
\usepackage{amsthm}
\usepackage{mathrsfs}
\usepackage{enumerate}
\newtheorem{theorem}{Theorem}
\theoremstyle{plain}

\newtheorem{corollary}{Corollary}

\newtheorem{definition}{Definition}

\newtheorem{lemma}{Lemma}
\newtheorem{notation}{Notation}

\numberwithin{equation}{section}
\newenvironment{example}{\begin{trivlist}\item[] {\bf Example.}\space\space\nobreak\noindent\ignorespaces}{\smallskip\end{trivlist}}
\newenvironment{remark}{\begin{trivlist}\item[] {\bf Remark.}\space\space\nobreak\noindent\ignorespaces}{\smallskip\end{trivlist}}

\input xy
\xyoption{all}
\newcommand{\zed}{\mathbb{Z}}

\newcommand{\projective}{\mathbb{P}}

\newcommand{\tensor}{\otimes}

\newcommand{\sheafo}{\mathscr{O}}

\newcommand{\morphismsend}{\text{End}}
\newcommand{\spec}{\text{Spec}}

\newcommand{\rank}{\text{rk}}
\newcommand{\characteristic}{\text{char}}
\newcommand{\coker}{\text{coker}}
\newcommand{\proj}{\text{Proj}}
\newcommand{\frobenius}{\text{Frob}}
\newcommand{\Span}{\text{span}}
\newcommand{\rk}{\text{rk}}
\newcommand{\im}{\text{im}}
\begin{document}
\title[A note on the A-numbers and P-ranks of Kummer covers]{A note on the A-numbers and P-ranks of Kummer covers}
\author{Otto Johnston}
\address[Otto Johnston]
{Universiteit van Amsterdam \newline%
\indent Amsterdam, The Netherlands}%
\email{o.r.johnston@student.uva.nl}%
\date{October 10, 2007}
\begin{abstract}
We study the $a$-numbers and $p$-ranks of Kummer covers of the projective line, and we give bounds for these numbers.
\end{abstract}
\maketitle
\section{Introduction}
\par For some special curves, explicit formulas exist for the $p$-rank in terms of $p$, the degree of $C$, and the degree of the ramification divisor.  One of the most famous of these formulas is due to Deuring and Shafarevich and dates back to the 1940s (see \cite{safarevic}).  However, as Crew pointed out much later in \cite{crew}, such a formula is impossible for Kummer covers since even for elliptic curves the $p$-rank can vary with the other numbers fixed.  The same argument works equally well for the $a$-number.
\par We will study the $a$-numbers and $p$-ranks of Kummer covers.   Our method uses \u{C}ech cohomology to produce a natural basis of $H^1(C, \sheafo_C)$, and we calculate the action of Frobenius on this basis.  Using this action, we give bounds for the $a$-number and the $p$-rank.  This extends recent results of Elkin who used a similar method for a more specialized class of Kummer cover (see \cite{elkin}).
\par As an application, we recover Ekedahl's bound $g(C) < \frac{p + 1}{2}$ for superspecial hyperelliptic curves (see \cite{ekedahl}).  We show that there are numbers less than this upper bound that do not occur as the genus of such a curve.
\section{Kummer Covers}
In this section, the main result is a decomposition theorem for the induced action of Frobenius on the first cohomology group of a Kummer cover.
\begin{definition}
An irreducible projective smooth curve $C$ over a field $k$ is a Kummer cover of degree $n$ if there exists a finite separable morphism $\psi: C \rightarrow \projective^1_k$ of degree $n$ such that $K(C)/K(\projective^1_k)$ is a Kummer extension.
\end{definition}
This definition automatically assumes that the characteristic of $k$ does not divide $n$ and that $k$ contains the $n$th roots of unity.  For example, hyperelliptic curves over algebraically closed fields of characteristic not equal to $2$ are the Kummer covers of degree $2$.  We will need the following algebraic fact.
\begin{lemma}
Let $R$ be a noetherian unique factorization domain, and let $R[\alpha]$ be the cyclic extension of $R$ defined by a root $\alpha$ of the irreducible polynomial $z^n - u\prod_{j = 1}^{n - 1} f_j^j$, where $f_j \in R$ is square-free and $u \in R^*$.  The integral closure of $R[\alpha]$ is generated as an $R$-module by
\[
\frac{\alpha^k}{\prod_{j = 1}^{n - 1} f_j^{\lfloor jk/n \rfloor}}, \;\;\;\; k = 0,1,\dots,n-1.
\]
\label{generators}
\end{lemma}
\begin{proof}
See \cite{EV}.
\end{proof}
\noindent We use the previous result to find an affine cover of our curve.
\begin{lemma}
Let $C$ be a Kummer cover of degee $n$ over a field $k_0$.  After a base extension $k/k_0$, we can find a generator $y$ for the cyclic extension $K(C)/K(\projective^1_k)$ such that $y^n= f = u\prod_{j = 1}^{n - 1} f_j^j$ for $u \in k^*$ and $f_j \in k[x]$ separable.  Then $C$ has an affine cover consisting of two parts $U^\prime = \emph{\spec}\;A$ and $V^\prime = \emph{\spec}\;B$, where
\[
A = \sum_{i = 0}^{n - 1} \frac{y^i}{\prod_{j = 1}^{n-1} f_j^{\lfloor ji/n \rfloor}} \cdot k[x],\;\;\;\;\;\; B = \sum_{i = 0}^{n - 1} \frac{y^i}{\prod_{j = 1}^{n-1} f_j^{\lfloor ji/n \rfloor}} \cdot \frac{1}{x^{m_i}} \cdot k[1/x],
\]
and
\[
m_i = \lceil i\deg(f)/n \rceil - \sum_{j = 1}^{n - 1} \deg(f_j) \lfloor ji/n \rfloor.
\]
\label{cover}
\end{lemma}
\begin{proof}
Let $K(\projective^1_{k_0}) = k_0(x)$.  Since $K(C)/k_0(x)$ is a Kummer extension, we can find a generator $\alpha$ such that $\alpha^n = q \in k_0[x]$ and $Z^n - q \in k_0(x)[Z]$ is irreducible.  We can also find a field extension $k/k_0$ such that all square-free factors of $q$ in $k[x]$ are separable.  Base extending $C$ and $\projective^1_{k_0}$ by $k$, we get a Kummer extension $K(C)/k(x)$ with a generator $y$ such that $y^n = f$, $f$ divides $q$, and $Z^n - f \in k(x)[Z]$ is irreducible.  We have the square-free factorization $f = u\prod_{j = 1}^{n - 1} f_j^j \in k[x]$, where each $f_j$ is separable since it divides $q$.
\par Write $\projective^1_k = \proj\;[t_0, t_1]$ and cover it by $U = \spec\;k[x]$ and $V = \spec\;k[1/x]$ with $x = t_1/t_0$.  Using the finite morphism $\psi: C \rightarrow \projective^1_k$, we form a cover of $C$ consisting of two affine open sets $U^\prime = \psi^{-1}(U)$ and $V^\prime = \psi^{-1}(V)$.
\par We know that $A = \Gamma(U^\prime, \sheafo_C)$ is the integral closure of $k[x]$ in $K(C)$ and $B = \Gamma(V^\prime, \sheafo_C)$ is the integral closure of $k[1/x]$ in $K(C)$ since $C$ is isomorphic to its normalization.  Lemma \ref{generators} immediately gives us the generators for $A$.  To find the generators for $B$, let $\alpha = (y/x^s)$ be the root of the irreducible polynomial $Z^n - u (\prod_{j = 1}^{n - 1} f_j)/x^{n s} \in k[1/x, Z]$, where $s = \lceil \deg(f)/n \rceil$.  We can use Lemma \ref{generators} to compute the integral closure of $k[1/x, \alpha]$.  Since the integral closure of $k[1/x]$ in $K(C)$ is the smallest integrally closed ring in $K(C)$ that contains $k[1/x]$, this computation is all we need.  Rearranging the basis elements for $B$ using elementary algebra gives us the desired form where $m_i = i\;s - \lfloor i\;s - i\;\deg(f)/n \rfloor - \sum_{j = 1}^{n - 1} \deg(f_j) \lfloor ji/n \rfloor$.  We get the definition of $m_i$ used above by the equality $i\;s - \lfloor i\;s - i\;\deg(f)/n \rfloor =  \lceil i\;\deg(f)/n \rceil$.
\end{proof}
\begin{lemma}
Using the same notation, let $C$ be a Kummer cover of degree $n$ over $k$.
\begin{enumerate}
\item $H^1(C, \sheafo_C) = \sum_{i = 1}^{n - 1} \sum_{t = 1}^{m_i - 1} \frac{y^i}{\prod_{j = 1}^{n-1} f_j^{\lfloor ji/n \rfloor}} \cdot \frac{1}{x^t} \cdot k$. \vspace{0.2in}
\item The genus of $C$ is $g(C) = \left( \sum_{i = 1}^{n - 1} m_i \right) - n + 1$.  Moreover,
\[
0 \leq g(C) \leq \frac{1}{2}(n - 1)(\deg(f) - 1).
\]
\item  Let $\emph{\characteristic}(k) = p$.  The induced Frobenius map $F$ on $H^1(C, \sheafo_C)$ is determined by
\[
\frac{y^i}{\prod_{j = 1}^{n-1} f_j^{\lfloor ji/n \rfloor}}\cdot \frac{1}{x^t} \longmapsto \sum_{w = 1}^{m_{(pi\bmod{n})} - 1} c_w \cdot \frac{y^{(pi\bmod{n})}}{\prod_{j = 1}^{n-1} f_j^{\lfloor j \cdot (pi\bmod{n})/n \rfloor}} \cdot \frac{1}{x^w},
\]
where $c_w$ is the coefficient of $x^{pt - w}$ in $f^{\lfloor pi/n \rfloor}/\prod_{j = 1}^{n - 1} f_j^{p \lfloor ji/n \rfloor - \lfloor j \cdot (pi\bmod{n})/n \rfloor}$.
\end{enumerate}
\label{Kummercovers}
\end{lemma}
\begin{proof}
(1)  Let $R = \Gamma(U^\prime \cap V^\prime, \sheafo_C)$ and note that $R$ is the integral closure of $k[x, 1/x]$ in $K(C)$.  Lemma \ref{generators} tells us that $R$ is generated as a $k[x, 1/x]$-module by the same set of generators that formed $A$ as a $k[x]$-module.  We form the standard \u{C}ech complex 
\[
\xymatrix{
A \oplus B \ar[r]^{\;\;d} & R \ar[r] & 0
}
\]
and the result is immediate after taking the quotient.
\par (2)  The genus formula is obvious from part (1).  The upper bound for $g(C)$ comes from considering $f$ to be square-free: it is clear we obtain the largest possible $m_i$ in this case, and hence the largest possible $g(C)$ for fixed $n$ and $\deg(f)$.  An obvious lower bound for the genus of a Kummer cover with a square-free $f$ is obtained by replacing $m_i$ with $i\cdot \deg(f)/n$, which gives us $\frac{1}{2}(n - 1)(\deg(f) - 2) \leq g(C)$.  The upper bound comes from the basic numerical fact that $\sum_{i = 1}^{n - 1} (m_i - (i)\deg(f)/n) \leq (n - 1)/2$, which is added to the formula for the lower bound.
\par (3)  We can determine the action of $F$ on $H^1(C, \sheafo_C)$ by the action of Frobenius on the \u{C}ech complex $A \oplus B \rightarrow R$.
Since $F$ is semi-linear on $k$, it is completely determined by its action on the basis vectors of $H^1(C, \sheafo_C)$.  To determine the action of $F$ on a basis vector, let $\frobenius$ denote the absolute Frobenius map on $C$ and look at the following commutative diagram.
\[
\xymatrix{
R \ar[d]_{\coker(d)} \ar[r]^{\frobenius} & R \ar[d]^{\coker(d)} \\
H^1(C, \sheafo_C) \ar[r]^F & H^1(C, \sheafo_C)
}
\]
We have already computed the basis vectors of $H^1(C, \sheafo_C)$ as the images of elements of $R$ under $\coker(d)$ of the form
\[
\frac{y^i}{\prod_{j = 1}^{n-1} f_j^{\lfloor ji/n \rfloor}} \cdot \frac{1}{x^t}.
\]
To compute the action of $F$ on a basis vector of $H^1(C, \sheafo_C)$, we will simply apply $\frobenius$ to the above term of $R$ and then apply $\coker(d)$.  Applying $\frobenius$ obviously gives us
\[
\frac{y^i}{\prod_{j = 1}^{n-1} f_j^{\lfloor ji/n \rfloor}} \cdot \frac{1}{x^t} \mapsto \left( \frac{y^i}{\prod_{j = 1}^{n-1} f_j^{\lfloor ji/n \rfloor}} \cdot \frac{1}{x^t} \right)^p,
\]
and we have that the image is equal to
\[
\frac{f^{\lfloor p i/n \rfloor}}{\prod_{j = 1}^{n-1} f_j^{p\lfloor ji/n \rfloor - \lfloor j (pi \bmod{n})/n \rfloor}} \cdot \frac{y^{(pi\bmod{n})}}{\prod_{j = 1}^{n-1} f_j^{\lfloor j \cdot (pi\bmod{n})/n \rfloor}} \cdot \frac{1}{x^{pt}} 
\]
by elementary algebra.  If we let $Q_i$ denote the leftmost term, we see that $Q_i \in k[x]$ since $j \lfloor p i/n \rfloor \geq p\lfloor ji/n \rfloor - \lfloor j (pi\bmod{n})/n \rfloor$.  To finish the calculation, we take the image of the above expression under $\coker(d)$, which is clearly $0$ if $m_{(pi\bmod{n})} \leq 1$.  If $m_{(pi\bmod{n})} > 1$, the image is the sum of the terms
\[
\left[ c_w \cdot \frac{y^{(pi\bmod{n})}}{\prod_{j = 1}^{n-1} f_j^{\lfloor j \cdot (pi\bmod{n})/n \rfloor}} \cdot \frac{1}{x^w} \right]
\]
for $w = 1, \dots, m_{(pi\bmod{n})} - 1$ and $c_w$ the coefficient of $x^{pt - w}$ as a term of $Q_i$.
\end{proof}
\begin{remark}
The bounds given for $g(C)$ in (2) are sharp.  The lower bound occurs for curves with affine equations $y^n = x^j$ for $j > 0$.  The upper bound occurs for  all Kummer covers with an affine equation of the form $y^n = f(x)$, where $f$ is separable and $\deg(f)$ is coprime to $n$. Also, we have seen that the computation of the \u{C}ech map involves the polynomial
\[
Q_i = f^{\lfloor pi/n \rfloor}/\prod_{j = 1}^{n - 1} f_j^{p \lfloor ji/n \rfloor - \lfloor j \cdot (pi\bmod{n})/n  \rfloor} \in k[x].
\label{Q}
\]
It is important to note that the exponents of the $1/f_j$ terms may be negative.
\end{remark}
\vspace{0.1in} 
\par We now turn our attention to the $a$-number and $p$-rank of a Kummer cover.  To define these numbers, we will need some facts about semi-linear maps.  Recall that a semi-linear map of a $k$ vector space $L: V  \rightarrow V$ is an additive map satisfying $L(\lambda x) = \theta(\lambda) L(x)$ for some $\theta \in \morphismsend(k)$.   For any semi-linear map $L$, the set $\ker(L)$ is a vector space over $k$ and $\im(L)$ is a vector space over $\theta(k)$.  Since it is more desirable to view the image of $L$  as a vector space over $k$, we define $\im_k(L) = \im(L) \tensor_{\theta(k)} k$.  
\par Many of the decomposition theorems from linear algebra carry over to semi-linear maps.  Recall that Rank-Nullity holds for $L$ in the sense that $\dim_k \ker(L) = r$ if and only if $\dim_k \im_k(L) = n - r$.  We also have that $\ker(L^m)$ stabilizes for some $m \geq 0$, where the smallest such $m$ is denoted by $i(L)$ and called the index of $L$.  Finally, the Range-Nullspace decomposition tells us that $L|_{\ker(L^m)}$ is nilpotent and $\dim_k \im_k(L^m) = \dim_k \im_k(L^{m + 1})$.  Of course, the semi-linear map we are interested in is $F$ acting on $H^1(C, \sheafo_C)$, where $\theta$ is $\lambda \mapsto \lambda^p$ on $k$.
\par From this point on, we assume that $\characteristic(k) = p > 0$.   The semi-simple rank of $F$ is $\rank(F) = \dim_k \im_k (F)$.  The $a$-number $a(C)$ of a curve $C$ over $k$ is $a(C) = \dim_k \ker(F)$.  Rank-Nullity gives us the relation $\rk(F) = g(C) - a(C)$.  The $p$-rank $f(C)$ of $C$ is $f(C) = \rk(F^m)$ for any $m \geq i(F)$.  This is well-defined because $\ker(F^m)$ stabilizes.  Moreover, it is easy to see that $i(F) \leq g(C)$, so we can always take $m$ to be $g(C)$ in the definition of $f(C)$.  The integers $\rk(F)$, $a(C)$, and $f(C)$ are all between $0$ and $g(C)$.  The curve $C$ is called \emph{superspecial} if $F = 0$.
\par The partition of $\zed/n\zed$ into subsets via the action of multiplication by $p$ plays an important role in our next result.  We fix the notation for this as follows.  
\begin{notation}
Let $S = \zed/n\zed - \{0\}$ and let $G$ be the cyclic group $\{p^q: q \geq 0\} \subset (\zed/n\zed)^*$.  Consider the group action of $G$ on $S$ given by $p^q \cdot s = p^q\;s \bmod n$.  Let $S/G$ be the set of distinct orbits of this action.
\label{notationone}
\end{notation}
\begin{theorem}
Using the same notation, let $C$ be a Kummer cover over $k$ of degree $n$.  Set $B_i = \emph{\Span}_k\{ \frac{y^i}{\prod_{j = 1}^{n-1} f_j^{\lfloor ji/n \rfloor}}\cdot \frac{1}{x^t}\}_{t = 1}^{m_i - 1}$.
\begin{enumerate}
\item $F^q(B_i) \subset B_{(i p^q\bmod n)}$ for $q > 0$. \vspace{0.1in}
\item $\emph{\rank}(F) = \sum_{i = 1}^{n - 1} \emph{\rank}(F|_{B_i})$. \vspace{0.1in}
\item $f(C) = \sum_{\Omega \in S/G} \emph{\rank}(F^m|_{B_i})$, where $m \geq i(F)$ and $i \in \Omega$ is any element.
\end{enumerate}
\label{decomposition}
\end{theorem}
\begin{proof}
Since $\sum_{i = 1}^{n - 1} B_i = H^1(C, \sheafo_C)$ by the first part of Lemma \ref{Kummercovers}, (2) and (3) follow from (1), so we prove (1).  Part (3) of Lemma \ref{Kummercovers} tells us that the action of $F$ takes $m_{i} - 1$ basis vectors and maps them to $m_{(pi\bmod{n})} - 1$ number of basis vectors.  Since multiplication by $p$ defines a bijection from $\zed/n\zed$ to itself, the $m_i - 1$ number of vectors are the only terms to be mapped to the $m_{(pi\bmod{n})} - 1$ number of vectors.  This proves $F(B_i) \subset B_{(p i\bmod n)}$.  Iterating $F$ finishes the proof.
\end{proof}
\vspace{0.1in}
\begin{example}
Let $C$ be the Kummer cover defined by $y^{11} = x^2 (x + 1)$ over a field of characteristic $13$ that contains the $11$th roots of unity.  We will show that $a(C) = 1$ and $f(C) = 0$ using the theorem.  The orbit of $1$ under the action of $G$ on $S$ is $\{1, 2, 4, 8, 5, 10, 9, 7, 3, 6\}$.  Thus, $S/G$ consists of the single orbit $S$.  Moreover, the set  $\{4, 5, 8, 9, 10\}$ consists of all values of $i < 11$ where $m_i > 1$; since $m_i = 2$ for these values, $g(C) = 5$.  On the other hand, since $S/G$ consists of a single orbit and there is some $j$ for which $m_j = 1$, the image of $F^q(B_i)$ passes through a zero dimensional $B_j$ for some iteration $q$ for any $i$.  Hence, $f(C) = 0$.  Since $m_7 = 1$ and $9$ maps to $7$ in one iteration, $\rank(F|_{B_9}) = 0$, so we can compute $\rank(F)$ by taking the sum of $\rank(F|_{B_i})$ for $i \in \{4, 5, 8, 10\}$.  To determine $\rank(F|_{B_i})$, all we need to know is if the coefficient $a_{i, 12}$ of $x^{12}$ in $Q_i$ is zero or not.   A simple computation reveals $a_{4, 12} = 4$, $a_{5, 12} = 5$, $a_{8, 12} = 10$, and $a_{10, 12} = 3$.  Thus, $\rank(F) = 4$ and $a(C) = 1$.  We see that $C$ is an example of a curve of genus $5$ with $a$-number $1$ and $p$-rank $0$.
\end{example}
\section*{Bounds for the invariants}
\par Using Theorem \ref{decomposition}, we can easily produce the following bounds.  The group action plays an important role in the calculation of the $p$-rank. \vspace{0.1in}
\begin{corollary}
Using the notation introduced in Lemma \ref{cover} and Notation \ref{notationone}, let $C$ be a Kummer cover of degree $n$ over $k$. \vspace{0.1in}
\begin{enumerate}
\item $a(C) \geq 1 - n + \sum_{i = 1}^{n - 1} \max\{1, m_i - m_{(pi\bmod{n})} + 1\}$. \vspace{0.1in}
\item $f(C) \leq \sum_{\Omega \in S/G} \min_{i \in \Omega}\{m_i - 1\}$. \vspace{0.1in}
\end{enumerate}
\label{lowerbound}
\end{corollary}
\begin{proof}
(1)  Using part (1) of Theorem \ref{decomposition} and Rank-Nullity, we have the bound $\rank(F|_{B_i}) \leq \min\{\dim_k B_i, \dim_k B_{pi \bmod n}\} = \min\{m_i - 1, m_{(pi\bmod{n})} - 1\}$.  We get the lower bound for $a(C)$ by subtracting the upper bound of $\rank(F)$ from $g(C)$ given in part (2) of Lemma~\ref{Kummercovers}.
\par (2)  Taking iterations in part (1) and using the Range-Nullspace decomposition, we have $\rank(F|_{B_i}^m) \leq \min\{\dim_k B_i, \dim_k B_{pi \bmod n}, \dim_k B_{p^2 i\bmod n}, \dots \} =  \min_{i \in \Omega}\{m_i - 1\}$ where $\Omega$ is the action of $G$ on $i$.
\end{proof}
\vspace{0.05in}
\begin{example}
The upper bounds are sharp.  For instance, take $C$ to be the curve $y^6 = x^3 + x^2 + 1$ over a field $k$ of characteristic $5$ that contains the $6$th roots of unity.  In this case, $G = \{1, 5\}$ and $S/G$ consists of the orbits $\{3\}$, $\{1, 5\}$, and $\{2, 4\}$.  Only $i$ in $\{3, 4, 5\}$ satisfies $m_i > 1$, where $m_3 = m_4 = 2$ and $m_5 = 3$.  From this information alone, we obtain the following: $g(C) = 4$, $f(C) \leq 1$, $\rk(F) \leq 1$, $a(C) \geq 3$, and $F = F|_{B_3}$.  The action of $F$ on $B_3$ is easy to determine: it is multiplication by the coefficient $a_4$ of $Q_3$, which is $1$.  Thus, our bounds are all met.  We see that $C$ is an example of a curve of genus $4$ with $a$-number $3$ and $p$-rank $g - 3$.
\end{example}
\begin{corollary}
Using the notation introduced in Lemma \ref{cover}, let $C$ be a Kummer cover of degree $n$ over $k$.
\[
a(C) \leq 1 - n + \sum_{i = 1}^{n - 1} \min\{m_i, \max\{m_i - q_i + v_i, 1 + v_i \} \}.
\]
where
\[
q_i = \lfloor (\deg(Q_i) + m_{(pi\bmod{n})} - 1)/p\rfloor\;\;\;\text{and}\;\;\;v_i = \lfloor \deg(Q_i)/p \rfloor. \vspace{0.1in}
\]
\label{upperbound}
\end{corollary}
\begin{proof}
Our task is to compute a lower bound for $\rank(F|_{B_i})$. The entries in $F|_{B_i}$ come from the coefficients of the polynomial $Q_i$ as described by Lemma \ref{Kummercovers}.  Let $c$ denote the leading coefficient of $Q_i$.  We will exploit the following fact: when $c$ is used in $F|_{B_i}$, we can use row-reduction to easily see that it must contribute $1$ to the rank (indeed, any coefficient of $Q_i$ can be used at most once on any given row and all entries below those coming from $c$ are zero).  This means we get a lower bound for $\rank(F|_{B_i})$ by counting the minimal number of rows where $c$ must occur; we compute this number as follows.  The integer $v_i$ is the largest possible row of $F|_{B_i}$ where $c$ may not occur since $p v_i \leq \deg(Q_i)$.  The largest row of $F|_{B_i}$ where $c$ must occur is $q_i$ since $p\;(q_i + 1) - m_{(pi\bmod{n})} + 1 > \deg(Q_i)$.  Thus, 
\[
\max\{0, \min\{q_i - v_i, m_i - 1 - v_i\}\} \leq \rank(F|_{B_i}).
\]
We conclude by taking the sum over $i$ of this lower bound for $\rank(F|_{B_i})$ and subtracting it from $g(C)$ as we did before.
\end{proof}
\noindent This lower bound can be made much stronger for superelliptic curves, see \cite{elkin}.
\section*{Hyperelliptic curves}
In this section, we look at hyperelliptic curves over an algebraically closed field $k$.  Since $a(C)$ and $f(C)$ are invariants under separable base extension, the assumption that $k$ is algebraically closed is no loss of generality for our purposes.  Hyperelliptic curves are Kummer covers in every characteristic except $2$, so we only need to extend our results to characteristic $2$.
\begin{lemma}
Let $C$ be a hyperelliptic curve of genus $g = g(C)$ over an algebraically closed field $k$ of characteristic $2$.  Assume that $C$ is ramified at infinity.
\begin{enumerate}
\item $C$ has an affine cover consisting of two parts $U^\prime = \emph{\spec}\;A$ and $V^\prime = \emph{\spec}\;B$, where
\vspace{0.05in}
\[
A = k[x, y]/(y^2 + Qy - P),
\]
\[
B = k[\frac{1}{x}, \frac{y}{x^{g + 1}}]/(\frac{y^2}{x^{2g + 2}} + \frac{Qy}{x^{2g + 2}} - \frac{P}{x^{2g + 2}}),
\]
\newline
and where $Q, P \in k[x]$ satisfy $\deg(Q) \leq g$, $\deg(P) = 2g + 1$, and $Q$ is coprime to $(Q^\prime)^2 P + (P^\prime)^2$. \vspace{0.2in}
\item $H^1(C, \sheafo_C) = \sum_{i = 1}^{g} k \cdot y/x^i$.  The induced action of Frobenius is given by
\[
y/x^i \mapsto \sum_{j = 1}^{g} c_{i,j} y/x^j,
\]
where $c_{i, j}$ is the coefficient of $x^{2i - j}$ as a term of $Q$. \vspace{0.2in}
\item If $f(C) = 0$, then $a(C) = \lfloor \frac{g + 1}{2} \rfloor$.
\end{enumerate}
\end{lemma}
\begin{proof} (1)  See Proposition 7.4.24 of \cite{qing}.
\par (2)  We have that $R = \Gamma(U^\prime \cap V^\prime, \sheafo_C) = k[x, 1/x, y]/(y^2 + Qy - P)$.  The result follows by forming the standard \u{C}ech complex $A \oplus B \rightarrow R$ and passing to the quotient.  As for the action induced by Frobenius, if we square $y/x^i$ in $H^1(C, \sheafo_C)$, we have the coset relation $[(Qy - P)/x^{2i}] = [Qy/x^{2i}] = \sum_{j = 1}^{g} c_{i,j} y/x^j$, where $c_{i, j}$ is the coefficient of $Q$ as stated.
\par (3)  View $F$ as the $g \times g$ matrix $(c_{i, j})$ and use the notation $F[i, j]$ to denote $c_{i, j}$.  Part (3) tells us that the $c_{i, j}$ are coefficients of the polynomial $Q = \sum a_i x^i$ of degree at most $g$.  Since $f(C) = 0$, we also know that $F$ is nilpotent.  We have that $(F^n)[g, g] = a_g^n$, which forces $a_g = 0$.  Using this, we continue our elimination: we have that $(F^n)[g - 1, g - 1] = a_{g - 1}^n$ forces $a_{g - 1} = 0$, $(F^n)[g - 2, g - 2] = a_{g - 2}^n$ forces $a_{g - 2} = 0$, and so on, until we have $(F^n)[1, 1] = a_{1}^n$, which forces $a_1 = 0$.  Hence, the only $Q$ that satisfies a nilpotent $F$ is the constant $Q = a_0$.  It must be non-zero because $0$ is not coprime to $(P^\prime)^2$.  If $g$ is even, $a_0$ appears on $g/2$ rows.  If $g$ is odd, $a_0$ appears on $(g - 1)/2$ rows.
\end{proof}
\noindent A much stronger version of (3) has been proved by G. van der Geer (see Lemma 11.1 of \cite{vdg}).
\par Ekedahl's bound $g(C) < (p + 1)/2$ for superspecial hyperelliptic curves is an immediate consequence of part (3) and Corollary \ref{upperbound} when we take $C \rightarrow \projective^1_k$ to be ramified over infinity.  It is well-known that this bound is sharp.  What we want to know is if all the numbers below Ekedahl's bound occur as the genus of some superspecial hyperelliptic curve in characteristic $p$.  For $g(C) = 2$ and $p > 3$, such curves exist by a result of Ibukiyama, Katsura, and Oort in \cite{ibu-kat-oort}.  The case $g(C) = 3$ and $p > 5$ follows from a result in \cite{brock}.  Despite these early successes, we will show that there are gaps for genus $4$ in the next example by showing that there is no superspecial hyperelliptic curve of genus $4$ in characteristic $11$.
\vspace{0.15in}
\begin{example}
Assume that $C$ is a superspecial hyperelliptic curve of genus $4$ over an algebraically closed field of characteristic $11$. Use a fractional linear transformation of $C$ to force $0$ and infinity to be ramification points.  Using Lemma \ref{cover}, $C$ has an affine equation of the form
\[
y^2 = f(x) = a_1 x + \dots + a_9 x^9,
\]
with $a_1 \neq 0$ and $a_9 \neq 0$.  Lemma $\ref{Kummercovers}$ tells us that $f(x)^5 = \sum b_i x^i$ has $b_j = 0$ for $j \in \{ 7, 8, 9, 10, 18, 19, 20, 21\}$.  Since $0 = b_7 = 10 a_1^3 a_2^2 + 5a_1^4 a_3$, we have the relation $a_3 = -2a_2^2/a_1$.  Likewise, $a_4 = -3a_2^3/5a_1^2$, $a_5 = -6a_2^4/5a_1^3$, and $a_6 = -8a_2^5/5a_1^4$.  For $b_{18}$, we have
\[
\frac{7 a_2^{13}}{a_1^8} + \frac{8 a_2^7 a_7}{a_1^3} + 8a_1^2 a_2 a_7^2 + \frac{4a_2^6 a_8}{a_1^2} + 9a_1^3 a_7 a_8 = 0.
\]
This breaks down into three possible statements that we enumerate and eliminate below.
\par I.  $a_8 = 0$ and $-7a_2^{13} - 8a_1^5 a_2^7 a_7 - 8a_1^{10} a_2 a_7^2 = 0$.  This gives us 
\[
0 = b_{19} = \frac{9 a_2^{14}}{a_1^9} + \frac{4a_2^8 a_7}{a_1^4} + 3a_1 a_2^2 a_7^2 + \frac{4a_2^6 a_9}{a_1^2} + 9a_1^3 a_7 a_9.
\]
Since $a_9 \neq 0$, we have two subcases from the condition above.  We enumerate and eliminate them.
\par I.a.  $a_7 = -4a_2^6/9a_1^5$ and $a_{2}^{14} = 0$.  This case is eliminated because $0 = b_{21} = 10a_1^4 a_9^2$, so either $a_1$ or $a_9$ is zero, which is not possible.
\par I.b.  $a_9 = -(9a_2^{14} + 4a_1^5 a_2^8 a_7 + 3a_1^{10} a_2^2 a_7^2)/a_1^7(4 a_2^6 + 9a_1^5 a_7)$.  Returning to $b_{18} = 0$, this forces $a_2 = 0$, which in turn forces $a_9 = 0$.
\par II.  $a_7 = -4a_2^6/9a_1^5$ and $a_2 = 0$.  We have that the relation $b_{19} = 0$ yields $a_8 = 0$, and then $b_{21} = 0$ forces either $a_1$ or $a_9$ to be zero.  
\par III.  $a_8 = -(7a_2^{13} + 8a_1^5 a_2^7 a_7 + 8a_1^{10}a_2 a_7^2)/(a_1^6 (4a_2^6 + 9a_1^5 a_7))$.  The relation $b_{19} = 0$, gives us the following subcases.
\par III.a.  $a_9 = 0$ and $-5a_2^{14} - 7a_1^{10}a_2^2 a_7^2 = 0$.  This case is impossible because $a_9 \neq 0$.
\par III.b.  $a_7 = -4a_2^6/9a_1^5$ and $-5a_2^{14} - 7a_1^{10} a_2^2 a_7^2 = 0$.  This case is impossible because this definition of $a_7$ conflicts with the definition of $a_8$ (it causes a division by $0$).
\par III.c.  $a_9 = -(5a_2^{14} + 7a_1^{10} a_2^2 a_7^2)/(a_1^7 (4a_2^6 + 9a_1^5 a_7))$.  We have that 
\[
0 = b_{20} = \frac{7a_1^{25} a_2^{15} + 8a_1^{30} a_2^9 a_7 + 8a_1^{35} a_2^3 a_7^2}{a_1^{35}}.
\]
On the other hand,
\[
0 = b_{21} = \frac{5a_1^{24} a_2^{16} + a_1^{29} a_2^{10} a_7 + a_1^{34} a_2^4 a_7^2}{a_1^{35}}.
\]
Combining the two yields $a_2 = 0$, which forces $a_9 = 0$.
\end{example}

\end{document}